\documentclass[a4paper, reqno]{amsart}
\usepackage{amssymb, amsmath, amsthm, chngcntr, enumerate, mathabx, mathrsfs, mathtools, tikz, url}
\usepackage[numbers]{natbib}

\usetikzlibrary{arrows,calc}

\numberwithin{equation}{section}
\newtheorem{thm}{Theorem}

\newtheorem{corollary}[thm]{Corollary}
\newtheorem{lemma}[thm]{Lemma}
\newtheorem{proposition}[thm]{Proposition}

\newtheorem{rmk}[thm]{Remark}

\newtheorem*{theorem*}{Theorem}

\newcommand{\C}{\mathbb{C}}
\newcommand{\R}{\mathbb{R}}
\newcommand{\Z}{\mathbb{Z}}
\newcommand{\N}{\mathbb{N}}
\newcommand{\T}{\mathbb{T}}
\newcommand{\D}{\mathbb{D}}
\newcommand{\E}{\mathbb{E}}

\title[A multiplier inclusion theorem on product domains ]{A multiplier inclusion theorem on product domains}

\author{Odysseas Bakas}
\address{Room 4606, James Clerk Maxwell Building, University of Edinburgh, Peter Guthrie Tait Road, Edinburgh, EH9 3FD.}
\email{o.bakas@sms.ed.ac.uk}

\date{}

\begin{document}

\begin{abstract}In this note it is shown that the class of all multipliers from the $d$-parameter Hardy space $H^1_{\mathrm{prod}} (\T^d)$ to $L^2 (\T^d)$ is properly contained in the class of all multipliers from $L \log^{d/2} L (\T^d)$ to $L^2(\T^d)$. 
\end{abstract}

\maketitle

\section{Introduction}\label{intro}
Let $d$ be a positive integer. If $X$ is a subspace of $L^1 (\T^d)$, then we denote by $\mathcal{M}_{X \rightarrow L^2 (\T^d)}$ the class of all multipliers from $X$ to $L^2 (\T^d)$, namely the class $\mathcal{M}_{X \rightarrow L^2 (\T^d)}$ consists of all functions $m: \Z^d \rightarrow \C$ such that for every $f \in X$ one has $\sum_{(k_1, \cdots, k_d)  \in \Z^d} |m( k_1, \cdots, k_d) \widehat{f} (k_1, \cdots, k_d)|^2 < \infty$. 

In \cite{second_paper}, it was shown that the class of all multipliers from the (real) Hardy space $H^1 (\T)$ to $L^2 (\T)$ is properly contained in the class of all multipliers from $L \log^{1/2} L (\T)$ to $L^2 (\T)$. Our goal in this note is to extend this result to the multi-parameter setting. First of all, note that if $H^1_{\mathrm{prod}} (\T^d)$ denotes the $d$-parameter (real) Hardy space over the $d$-torus, then $ L \log^d L (\T^d) \subset  H^1_{\mathrm{prod}} (\T^d) $ and hence, one automatically has  $\mathcal{M}_{H^1_{\mathrm{prod}} (\T^d)  \rightarrow L^2 (\T^d)} \subset \mathcal{M}_{L \log^d L (\T^d) \rightarrow L^2 (\T^d)}$. On the other hand, by adapting the argument given in \cite{second_paper} to the multi-parameter case, one deduces that the best we can expect is that $\mathcal{M}_{H^1_{\mathrm{prod}} (\T^d) \rightarrow L^2 (\T^d)} \subset \mathcal{M}_{L \log^{d/2} L (\T^d) \rightarrow L^2 (\T^d)}$. In this note we prove that this is indeed the case, namely we strengthen the trivial exponent $r=d$ in $L \log^r L (\T^d)$ to the optimal one, $r=d/2$. In particular, our main result in this note is the following theorem.
\begin{thm}\label{mult_incl}
One has the inclusion
\begin{equation}\label{inclusion}
 \mathcal{M}_{H^1_{\mathrm{prod}} (\T^d ) \rightarrow L^2 (\T^d)} \subset \mathcal{M}_{L \log^{d/2} L (\T^d) \rightarrow L^2 (\T^d)}.
\end{equation}
Moreover, the above inclusion is proper and it is sharp, in the sense that the exponent $r=d/2$ in $L \log^{d/2} L (\T^d)$ cannot be improved.
\end{thm}

The multiplier inclusion (\ref{inclusion}) is obtained by a series of reductions. First, arguing as in \cite{second_paper} and by using D. Oberlin's characterisation of the class $\mathcal{M}_{H^1_{\mathrm{prod}} (\T^d) \rightarrow L^2 (\T^d)}$ given in \cite{Oberlin}, it follows that the proof of (\ref{inclusion}) is reduced to showing the following higher-dimensional version of an inequality due to Zygmund (see Theorem 7.6 in Chapter XII of \cite{Zygmund_book}), a result of independent interest. To state this version of Zygmund's inequality on $\T^d$, let $\mathcal{J}$ denote the set of all ``intervals'' of integers of the form $\pm \{ 2^n -1, \cdots, 2^{n+1}-2 \}$, $n \in \N_0$, in other words, $\mathcal{J}$ consists of all the sets in $\Z$ of the form $\{ 2^k -1, \cdots, 2^{k+1}-2 \}$, $k \in \N_0$ and $\{ -2^{l+1} +2, \cdots, -2^l +1 \}$, $l \in \N_0$.

\begin{proposition}\label{2d_Zygmund}
Let $\mathcal{J}$ be as above.

If $E \subset \Z^d$ is a non-empty set satisfying the condition
\begin{equation}\label{rectangle}
 D_E =\sup_{I_1 , \cdots, I_d \in \mathcal{J}} \# \big\{ E \cap (I_1 \times \cdots \times I_d) \big\} < \infty,
\end{equation}
then there exists a positive constant $A_{D_E}$, depending only on $D_E$, such that
\begin{equation}\label{non-product}
\Big(  \sum_{(k_1, \cdots, k_d)\in E} |\widehat{f} (k_1, \cdots, k_d)|^2 \Big)^{1/2} \leq A_{D_E} [ 1+ \int_{\T^d } |f| \log^{d/2}(1+|f|)]  .
\end{equation}

\end{proposition}

In turn, (\ref{non-product}) will be a corollary of a higher-dimensional extension of a result due to Seeger and Trebels \cite{Seeger_Trebels} concerning sharp bounds of sums involving ``smooth'' Littlewood-Paley projections on $\T^d$. To state this result, fix a Schwartz function $\eta$ supported in $(-2,2)$ such that $\eta|_{[-1,1]} \equiv 1$ and consider $\phi (\xi) = \eta(\xi) - \eta (2 \xi)$. For $k \in \N$, set $\phi_k (\xi) = \phi(2^{-k} \xi)$ and for $k=0$, set $\phi_0 = \eta$. One can easily see that $\sum_{k \in \N_0} \phi_k (\xi) = 1$ for every $\xi \in \R$.  Then, for $k \in \N_0$, the corresponding ``smooth'' Littlewood-Paley projection in the periodic setting is defined by  
$$ \widetilde{\Delta}_k (f) (x) = \sum_{r \in \Z} \phi_k (r) \widehat{f} (r) e^{i2 \pi r x} $$
for any, say, trigonometric polynomial $f$ on $\T$. On the $d$-torus we put
\begin{align*} 
\widetilde{\Delta}_{k_1, \cdots, k_d }(f) (x_1, \cdots, x_d) &= \widetilde{\Delta}_{k_1} \otimes \cdots \otimes \widetilde{\Delta}_{k_d} (f) (x_1, \cdots, x_d) \\
&=  \sum_{r_1 ,  \cdots, r_d \in \Z} \phi_{k_1} (r_1) \cdots \phi_{k_d} (r_d) \widehat{f} (r_1, \cdots, r_d) e^{i2 \pi (r_1 x_1 + \cdots + r_d x_d)} 
\end{align*}
initially defined over trigonometric polynomials $f$ on $\T^d$.
Then, Proposition \ref{2d_Zygmund} is a consequence of the following result.

\begin{proposition}\label{LP_product} There exists a constant $C_d>0$, depending only on the dimension $d$ and our choice of $\phi$, such that the following inequality holds
\begin{equation}\label{variant_LP}
\| f \|_{L^p (\T^d)} \leq C_d p^{d/2}  \Big(  \sum_{k_1, \cdots, k_d \in \N_0} \| \widetilde{\Delta}_{k_1, \cdots, k_d}(f) \|^2_{L^{\infty}(\T^d)}  \Big)^{1/2}
\end{equation}
for every trigonometric polynomial $f$ on $\T^d$ and for each $p>2$. 
\end{proposition}

The proof of Proposition \ref{LP_product} is an adaptation of the work of Seeger and Trebels \cite{Seeger_Trebels} to the higher-dimensional setting combined with a well-known inequality on multiple martingales, see section \ref{martingale_ineq}. At this point, it should be mentioned that, in fact, we expect that
$$ \| f \|_{L^p (\T^d)} \lesssim_d p^{d/2}  \Big\| \Big(  \sum_{k_1, \cdots, k_d \in \N_0} | \widetilde{\Delta}_{k_1, \cdots, k_d}(f) |^2   \Big)^{1/2} \Big\|_{L^p (\T^d)}$$
 which, of course,  implies (\ref{variant_LP}). However, as our primary goal is to establish Theorem \ref{mult_incl} and since (\ref{variant_LP}) is enough for that purpose, we shall not pursue this in the present note. 

The paper is organised as follows.
In section \ref{notation} we give some notation and background and in section \ref{2,1} we show how the proof of our multiplier inclusion theorem follows from Proposition \ref{2d_Zygmund}. In section \ref{3,2}, we prove that Proposition \ref{LP_product} implies Proposition \ref{2d_Zygmund} and then, in section \ref{proof} we give a proof of Proposition \ref{LP_product}. In the last section we briefly present some further applications of our work.

\subsection*{Acknowledgement} The author would like to thank his PhD supervisor Professor Jim Wright for his guidance on this work and for his useful comments that improved the presentation of this paper.

\section{Notation and background}\label{notation}

We denote by $\Z$ the set of integers, by $\N$ the set of positive integers, and by $\N_0$  the set of non-negative integers.

The cardinality of a finite set $A$ is denoted by $\#\{A\}$.

If $X$ and $Y$ are positive quantities such that $X \leq C Y$, where $C>0$ is a constant, then we write $X \lesssim Y$. To specify the  dependence of this constant on some additional parameters $\alpha_1, \cdots, \alpha_n$ we write $ X \lesssim_{\alpha_1, \cdots, \alpha_n} Y$. If $X \lesssim Y$ and $Y \lesssim X$, we write $X \sim Y$.

In this note, we identify $\T$ with $[0,1)$ in the standard way.
\subsection{Product Hardy spaces and the class $\mathcal{M}_{H^1_{\mathrm{prod}} (\T^d ) \rightarrow L^2 (\T^d)}$}
For $0<r<1$, let $P_r$ denote the Poisson kernel on $\T$ given by $ P_r (x) = (1-r^2 )/(1-2r \cos x + r^2)$, $x \in \T$. For $x \in \T$, let $\Gamma(x) = \{ z \in \D: |z - e^{i 2 \pi x}| \leq 2 (1-|z|) \}$, where $\D$ denotes the unit disc in the complex plane. Then, the $d$-parameter (real) Hardy space $H^1_{\mathrm{prod}} (\T^d)$ consists of all integrable functions $f $ on the $d$-torus such that $f^{\ast} \in L^1 (\T^d)$, where for $(x_1, \cdots, x_d) \in \T^d$ one has
$$ f^{\ast} (x_1, \cdots, x_d) = \sup_{ z_1 \in \Gamma (x_1) , \cdots, z_d \in \Gamma (x_d)} |f \ast ( P_{r_1} \otimes \cdots \otimes P_{r_d} ) (z_1, \cdots, z_d)| .$$

It follows by the work of D. Oberlin \cite{Oberlin} that $m:\Z^d \rightarrow \C$ belongs to the class $ \mathcal{M}_{H^1_{\mathrm{prod}} (\T^d) \rightarrow L^2 (\T^d)}$ if and only if, 
\begin{equation}\label{Oberlin_cond}
\sup_{N_1, \cdots, N_d \in \N_0} \sum_{N_1 \leq |k_1| \leq 2 N_1} \cdots \sum_{N_d \leq |k_d| \leq 2 N_d} |m(k_1, \cdots, k_d)|^2 < \infty . 
\end{equation}
\subsection{Dyadic square functions}\label{martingale_ineq}
If $f\in L^1(\T)$ and $m\in \N_0$, then the $m$-th conditional expectation of $f$ is given by
$$ \E_m (f) (x) = 2^m \int_I f(x') dx',  $$
where $I $ is the unique dyadic interval in $\T$ of the form $I  = [s  2^{-m}, (s+1)2^{-m})$, $s=0,1,\cdots, 2^m-1$ such that $x \in I $. 

For $m \in \N$, let $\D_m = \E_m - \E_{m-1} $ denote the martingale differences acting on functions defined on $\T$. For $m=0$, we set $\D_0=\E_0$. 

For a given $d$-tuple $(m_1, \cdots, m_d )$ of non-negative integers, we define
$$ \E_{m_1, \cdots, m_d} = \E_{m_1} \otimes \cdots \otimes \E_{m_d} $$ 
and 
$$\D_{m_1, \cdots ,  m_d} = \D_{m_1} \otimes \cdots \otimes \D_{m_d} = (\E_{m_1} - \E_{m_1 - 1}) \otimes \cdots \otimes (\E_{m_d} - \E_{m_d - 1})$$
 to be the corresponding operators acting on functions on the $d$-torus. 

In \cite{CWW}, Chang, Wilson, and Wolff obtained the`` good-$\lambda$'' inequality
\begin{align*}
\big| \big\{ x \in \T: \sup_{m \in \N_0} |\E_m f(x)| > 2\lambda \big\} \cap\big\{ x \in \T: \big( \sum_{m \in \N_0} |\D_m (f) (x) |^2 \big)^{1/2} < \epsilon \lambda \big\} \big|& \leq \\
 C_0 \exp[-{ \frac{(1-\epsilon)^2} {2 \epsilon^2} } ] \big| \big\{ x \in \T : \sup_{ m \in \N_0} |\E_m f(x)| > \lambda \big\} \big|&,
\end{align*}
which holds for all $\lambda>0$ and $\epsilon>0$, where $C_0 >0$ is an absolute constant. In particular, this estimate implies that there exists a constant $C>0$ such that
\begin{equation}\label{CWW_ineq}
\| f \|_{L^p (\T)} \leq C p^{1/2} \big\|  \big( \sum_{m \in \N_0} |\D_m (f) (x) |^2 \big)^{1/2} \big\|_{L^p (\T)}
\end{equation}
for all $p>2$.  By using (\ref{CWW_ineq}), Chang, Wilson, and Wolff obtained in \cite{CWW} an inequality analogous to (\ref{CWW_ineq}) involving Lusin area integrals. In \cite{Pipher}, Pipher extended (\ref{CWW_ineq}) and its analogous version on Lusin area integrals to the two-parameter setting and in \cite{Fefferman_Pipher}, R. Fefferman and Pipher extended the aforementioned inequality of Chang, Wilson, and Wolff involving Lusin area integrals to $\ell^2$-valued functions. The argument of R. Fefferman and Pipher can be easily adapted to obtain an $\ell^2$-valued extension of (\ref{CWW_ineq}), see \cite{Demeter_DiPlinio}. By using this $\ell^2$-valued extension of (\ref{CWW_ineq}) together with induction on $d$, one deduces that there exists a constant $C_d >0$, depending only on the dimension $d \in \N$, such that
\begin{equation}\label{Pipher}
  \| f \|_{L^p (\T^d)} \leq C_d p^{d/2} \Big\| \Big(  \sum_{m_1, \cdots , m_d \in \N_0} |\D_{m_1, \cdots , m_d} (f) |^2 \Big)^{1/2} \Big\|_{L^p (\T^d)}  
\end{equation}
 for every $p>2$, see also, e.g., \cite[Proposition 4.5]{Demeter_DiPlinio} and \cite{Bilyk}.
\subsection{Thin sets in Harmonic Analysis}\label{thin_sets} 
Let $G$ be a compact abelian group and let $\Lambda$ be a non-empty set in its dual $\widehat{G}$. In this note, we shall only consider the case $G = \T^d$, $d \in \N$. A trigonometric polynomial $f$ on $G$ whose spectrum lies in $\Lambda$ is said to be a $\Lambda$-polynomial.
 
Let $p >2$. We say that $\Lambda \subset \widehat{G}$ is a $\Lambda(p) $ set if there exists a constant $A(\Lambda, p) >0$ such that
$$  \| f \|_{L^p (G)} \leq A(p,\Lambda) \| f \|_{L^2 (G)}$$
for every $\Lambda$-polynomial $f$. The smallest constant $A(p,\Lambda)$ such that the above inequality holds is called the $\Lambda(p) $ constant of $\Lambda$.

A set $\Lambda \subset \widehat{G}$ is called Sidon if there is a constant $S_{\Lambda}>0$ such that
\begin{equation}\label{Sidon}
 \sum_{\gamma \in \Lambda} |\widehat{f} (\gamma)| \leq S_{\Lambda} \| f\|_{L^{\infty}(G)}
\end{equation}
for every $\Lambda$-polynomial. It follows by the work of Rudin \cite{Rudin_gaps} and Pisier \cite{Pisier_Sidon} that a spectral set $\Lambda $ is Sidon if and only if, it is a $\Lambda(p) $ set for any $p>2$ and its $\Lambda(p)$ constant grows like $p^{1/2}$ as $p \rightarrow \infty$.

Let $q \geq 1$. A set $\Lambda \subset \widehat{G}$ is said to be $q$-Rider if there is a constant $R_{\Lambda,q}>0$ such that
\begin{equation}\label{Rider}
  \big( \sum_{\gamma \in \Lambda} |\widehat{f} (\gamma)|^q \big)^{1/q} \leq R_{\Lambda,q} [| f |] 
\end{equation}
for every $\Lambda$-polynomial.  Here, we use the notation $[| f |] = \mathbb{E}\big[   \Big\| \sum_{\gamma \in \widehat{G}} r_{\gamma}  \widehat{f} (\gamma) \gamma \Big\|_{L^{\infty}(G)} \big]$, where $(r_{\gamma})_{\gamma}$ denotes the set of Rademacher functions. 

It is well-known that if $\Lambda$ is a $\Lambda (p)$ set for all $p>2$ with $\Lambda(p)$ constant growing like $p^{k/2}$, $k \in \N$, then $\Lambda$ is a $q$-Rider set with $q=2k/(k+1)$, see \cite[Th\'eor\`eme 6.3]{Pisier_aleatoires}.

\section{Proposition \ref{2d_Zygmund} implies Theorem \ref{mult_incl}}\label{2,1} 
To prove that Proposition \ref{2d_Zygmund} implies Theorem \ref{mult_incl}, we adapt the argument given in \cite{second_paper} to the multi-parameter setting by using the characterisation of $\mathcal{M}_{H^1_{\mathrm{prod}} (\T^d ) \rightarrow L^2 (\T^d)}$.  To be more specific, assume that Proposition \ref{2d_Zygmund} holds and take an arbitrary $m $ in the class $ \mathcal{M}_{H^1_{\mathrm{prod}}(\T^d) \rightarrow L^2 (\T^d)} $. Then, by definition, we need to show that for every $f \in L \log^{d/2} L(\T^d)$ one has 
$$\sum_{(k_1, \cdots, k_d)\in \Z^d} |m(k_1, \cdots, k_d) \widehat{f} (k_1, \cdots, k_d)|^2 <\infty. $$ 
Towards this aim, fix an $f \in L \log^{d/2} L (\T^d) $ and note that the sum
$$\sum_{(k_1, \cdots , k_d) \in \Z^d} |m(k_1, \cdots, k_d) \widehat{f} (k_1,\cdots, k_d) |^2  $$ 
is bounded by
$$ \sum_{I_1, \cdots, I_d \in \mathcal{J}} \max_{(k_1, \cdots, k_d) \in I_1 \times \cdots \times I_d} |\widehat{f} (k_1, \cdots, k_d)|^2  \big( \sum_{k_1 \in I_1} \cdots \sum_{k_d \in I_d} |m(k_1, \cdots, k_d) |^2 \big) , $$
where $\mathcal{J}$ is as in the introduction and the statement of Proposition \ref{2d_Zygmund}.
Hence, by (\ref{Oberlin_cond}), it follows that
$$\sum_{(k_1, \cdots, k_d)\in \Z^d} |m(k_1, \cdots, k_d) \widehat{f} (k_1, \cdots, k_d)|^2   \lesssim_m \sum_{(\widetilde{k}_1, \cdots , \widetilde{k}_d) \in E_f}  |\widehat{f} (\widetilde{k}_1, \cdots , \widetilde{k}_d) |^2, $$
where $E_f $ is a set in $\Z^d$ defined as follows. Given $I_1, \cdots, I_d \in \mathcal{J}$, choose $(\widetilde{k}_1 , \cdots, \widetilde{k}_d) $ in $ I_1 \times \cdots \times I_d$  so that

$$ |\widehat{f} (\widetilde{k}_1, \cdots, \widetilde{k}_d)| =  \max_{ (k_1, \cdots, k_d) \in I_1 \times \cdots \times I_d } |\widehat{f} (k_1, \cdots, k_d)|.$$
Then, having chosen a set of $d$-tuples $(\widetilde{k}_1 , \cdots, \widetilde{k}_d) $ as above, we define 
$$E_f = \{ (\widetilde{k}_1 , \cdots, \widetilde{k}_d) \in \Z^d: \ \mathrm{for} \ I_1 ,\cdots,  I_d \in \mathcal{J}, \ (\widetilde{k}_1, \cdots , \widetilde{k}_d)\in I_1 \times \cdots \times I_d\ \mathrm{being}\ \mathrm{as}\ \mathrm{above} \} .$$ 
Notice that as the choice of $d$-tuples $(\widetilde{k}_1, \cdots , \widetilde{k}_d)$ is not necessarily unique, there might be several choices of sets $E_f$. We just choose one of them to write
$$  \sum_{I_1, \cdots, I_d \in \mathcal{J}} \max_{(k_1, \cdots, k_d) \in I_1 \times \cdots \times I_d} |\widehat{f} (k_1, \cdots, k_d)|^2 =   \sum_{(\widetilde{k}_1, \cdots , \widetilde{k}_d) \in E_f}  |\widehat{f} (\widetilde{k}_1, \cdots , \widetilde{k}_d) |^2 .$$
Note that any such set $E_f$ satisfies condition (\ref{rectangle}) in Theorem \ref{2d_Zygmund} with $D_{E_f} = 1$. Therefore, as $f \in L\log^{d/2} L (\T^d)$, it follows by (\ref{non-product}) that
$$\sum_{(k_1, \cdots, k_d)\in \Z^d} |m(k_1, \cdots, k_d) \widehat{f} (k_1, \cdots, k_d)|^2 <\infty ,$$
as desired.

\subsection{Sharpness of (\ref{inclusion})}
We remark that, in fact, the above argument shows that if $m \in \mathcal{M}_{H^1_{\mathrm{prod}}(\T^d) \rightarrow L^2 (\T^d)}$, then there is a constant $C_m>0$, depending only on $m$, such that
$$ \Big( \sum_{(k_1, \cdots, k_d) \in \Z^d} |m(k_1, \cdots , k_d) \widehat{f} (k_1, \cdots, k_d)|^2 \Big)^{1/2} \leq C_m \big[ 1 + \int_{\T^d} |f| \log^{d/2} (1+|f|) \big]. 
 $$

To see that the exponent $r=d/2$ in $L \log^{d/2} L (\T^n)$ in (\ref{inclusion}) cannot be improved, we argue as in \cite{second_paper}. More specifically,  assume that for some $r>0$ every multiplier from $H^1_{\mathrm{prod}}(\T^d)$ to $L^2 (\T^d)$ is a multiplier from $L \log^r L (\T^d)$ to $L^2 (\T^d)$. We shall prove that $r \geq  d/2$. To this end, for a large positive integer $N$, take $f$ to be a trigonometric polynomial on $\T^d$ given by $f = V_{2^N} \otimes \cdots \otimes V_{2^N}$, where
 $V_{2^N} = 2 K_{2^{N+1} - 1} - K_{2^N - 1} $ denotes the de la Vall\'{e}e Poussin kernel of order $2^N$ and $K_n$ is the Fej\'er kernel on $\T$ of order $n \in \N$. Since $\| K_n \|_{L^1 (\T)} = 1$ and $\| K_n \|_{L^{\infty} (\T)}  \lesssim n$, we deduce that
  $$ \int_{\T^d} |f (x_1 ,\cdots, x_d) | \log^r (1+|f (x_1, \cdots, x_d)|) d x_1 \cdots d x_d \lesssim_{r,d} N^r.$$
So, if we take $M = (m(k_1, \cdots, k_d))_{k_1, \cdots, k_d \in \Z}$ with $m(k_1, \cdots, k_d)=1/\sqrt{k_1 \cdots k_d}$ for $k_1 > 0, \cdots, k_d >0$ and $m(k_1, \cdots, k_d)=0$ otherwise, namely when at least one of the coordinates is less or equal than $0$, then $M \in \mathcal{M}_{H^1_{\mathrm{prod} }(\T^d) \rightarrow L^2 (\T^d)}$ and hence,
$$ \Big( \sum_{(k_1, \cdots, k_d) \in \Z^d} |m (k_1, \cdots , k_d) \widehat{f}(k_1, \cdots, k_d) |^2 \Big)^{1/2} \lesssim_{r,d} N^r .$$
Since
\begin{align*}
 \Big( \sum_{(k_1, \cdots, k_d) \in \Z^d} |m (k_1, \cdots, k_d) \widehat{f}(k_1, \cdots, k_d) |^2 \Big)^{1/2} &\geq ( \sum_{1 \leq k_1, \cdots, k_d \leq 2^N} \frac{1}{k_1  \cdots k_d} \Big)^{1/2} \\
 &= \prod_{i=1}^d \Big(  \sum_{1 \leq k_i \leq 2^N} \frac{1}{k_i } \Big)^{1/2}   \\
 &\sim N^{d/2} , 
 \end{align*}
we see that, by choosing $N$ to be large enough, we must have $r \geq d/2$.

\begin{rmk}
A similar argument shows that the Orlicz space $L \log^{d/2} L (\T^d)$ in $(\ref{non-product})$ cannot be improved. Indeed, if $E$ is a set satisfying $(\ref{rectangle})$, then by making use of the argument presented above, we see that the exponent $r=d/2 $ in $L \log^{d/2} L (\T^d)$ in the right-hand side of higher-dimensional Zygmund's inequality $(\ref{non-product})$ is sharp.
\end{rmk}

To show that the inclusion $(\ref{inclusion})$ is proper, take $\Lambda $ to be a Sidon set in $ \Z $ that cannot be written as a finite union of lacunary sequences, see \cite[Remark 2.5(3)]{Rudin_gaps}. 
Then $ M = \chi_{\Lambda \times  \cdots \times \Lambda} $ belongs to the class $\mathcal{M}_{L \log^{d/2} L (\T^d) \rightarrow L^2 (\T^d)}$, see, e.g., \cite[Proposition 4]{second_paper}. However, it can be easily checked that $M = \chi_{\Lambda \times  \cdots \times \Lambda} $ does not satisfy (\ref{Oberlin_cond})
and hence, we deduce that $\chi_{\Lambda \times \cdots \times \Lambda} \in \mathcal{M}_{L \log^{d/2} L (\T^d) \rightarrow L^2 (\T^d)} \setminus \mathcal{M}_{H^1_{\mathrm{prod}} (\T^d) \rightarrow L^2 (\T^d)}$.

\section{Proposition \ref{LP_product} implies Proposition  \ref{2d_Zygmund}}\label{3,2}
Our goal in this section is to prove that Proposition \ref{LP_product} implies Proposition \ref{2d_Zygmund}. Towards this aim, take $E \subset \Z^d$ to be a set satisfying the assumption of Proposition \ref{2d_Zygmund}, i.e. condition (\ref{rectangle}). Assume first that $E$ satisfies (\ref{rectangle}) with $D_E =1$. By duality, to prove (\ref{non-product}), it suffices to show that $E$ is a $\Lambda(p)$ set in $\Z^d$ for every $p>2$ with $\Lambda(p)$ constant growing like $p^{d/2}$ as $p \rightarrow \infty$. In other words, it is enough to show that for every $E$-polynomial $f$  one has for every $p>2$,
\begin{equation}\label{dual}
\| f \|_{L^p (\T^d) } \leq A_E p^{d/2} \| f \|_{L^2 (\T^d)},
\end{equation}
where $A_E$ is an absolute constant, independent of $p$ and $f$. As we will see momentarily, if $D_E=1$, then, in fact, $A_E$ depends only on $d$ and in particular, it can be taken to be independent of $E$. 

To prove (\ref{dual}), fix an $E$-polynomial $f$ and note that for every $(k_1, \cdots, k_d)\in \N^d_0$ one has by the triangle inequality
\begin{align*}
\| \widetilde{\Delta}_{k_1, \cdots, k_d}(f) \|_{L^{\infty} (\T^d)} &\leq \sum_{(r_1, \cdots, r_d) \in E \cap  ( I_{k_1} \times \cdots \times I_{k_d} )  } |\phi_{k_1} (r_1) \cdots \phi_{k_d} (r_d) \widehat{f} (r_1, \cdots , r_d) | \\
&\lesssim_{\phi} \sum_{(r_1, \cdots, r_d) \in E \cap ( I_{k_1} \times \cdots \times I_{k_d}  ) }  | \widehat{f} (r_1, \cdots , r_d) | ,
\end{align*}
where $I_{k_l}$ denotes the set $\Z \cap \{ (-2^{k_l +1}, -2^{k_l-1}] \cup [2^{k_l-1},2^{k_l+1}) \} $, $l=1, \cdots, d$. Observe that, thanks to condition (\ref{rectangle}) for $D_E=1$, the sum 
$$ \sum_{(r_1, \cdots, r_d) \in E \cap ( I_{k_1} \times \cdots \times I_{k_d}  ) }  | \widehat{f} (r_1, \cdots , r_d) | $$ consists of at most $6^d$ terms and hence,
$$ \| \widetilde{\Delta}_{k_1, \cdots, k_d}(f) \|^2_{L^{\infty} (\T^d)} \lesssim_{d,\phi} \sum_{(r_1, \cdots, r_d) \in E \cap ( I_{k_1} \times \cdots \times I_{k_d}  ) }  | \widehat{f} (r_1, \cdots , r_d) |^2 $$
and we thus deduce that
\begin{equation}\label{bound_square} \Big(  \sum_{k_1, \cdots, k_d \in \N_0} \| \widetilde{\Delta}_{k_1, \cdots, k_d}(f) \|^2_{L^{\infty}(\T^d)}  \Big)^{1/2} \lesssim_{d, \phi}  \Big(  \sum_{(r_1, \cdots, r_d) \in E  }  | \widehat{f} (r_1, \cdots , r_d) |^2 \Big)^{1/2}.
\end{equation}
Observe that the quantity on the right-hand side of the last inequality equals to $\| f \|_{L^2 (\T^n)}$, as $\mathrm{supp}(\widehat{f}) \subset E$. Hence, (\ref{dual}) follows from (\ref{variant_LP}) and (\ref{bound_square}) in the case where $D_E =1$. Moreover, note that, in the case where $D_E=1$, the implied constant in (\ref{bound_square}) depends only on the dimension $d$ and on our choice of $\phi$ and, in particular, it is independent of $E$.

In the case where $D_E >1$, write $f =\sum_{i=1}^{D_E} f_i$, where $f_i$ are trigonometric polynomials on $\T^d$ such that $\mathrm{supp}(\widehat{f_i}) \subset E_i$, where $E= \cup_{i=1}^{D_E} E_i$ and $D_{E_i} =1$. Then, by using the triangle inequality and the previous step we have
$$ \| f \|_{L^p (\T^d)} \leq \sum_{i=1}^{D_E} \| f_i \|_{L^p (\T^d)} \leq A p^{d/2} \sum_{i=1}^{D_E} \| f_i \|_{L^2 (\T^d)} \leq A D_E p^{d/2} \| f \|_{L^2 (\T^d)},  $$
since, by our construction and the $L^2$-theory, $\| f_i \|_{L^2 (\T^d)} \leq \| f \|_{L^2 (\T^d)}$ for all $i=1, \cdots , D_E$.

\section{Proof of Proposition \ref{LP_product}}\label{proof}
To prove Proposition \ref{LP_product}, note that, as $p>2$, it follows by Minkowski's inequality that
$$  \Big\| \Big(  \sum_{m_1, \cdots, m_d \in \N_0} |\D_{m_1, \cdots, m_d} (f) |^2 \Big)^{1/2} \Big\|_{L^p (\T^d)} \leq   \Big(  \sum_{m_1, \cdots, m_d \in \N_0} \|\D_{m_1, \cdots, m_d} (f) \|^2_{L^p (\T^d)} \Big)^{1/2}.  $$
Moreover, since one trivially has
$$\Big(  \sum_{m_1, \cdots, m_d \in \N_0} \|\D_{m_1, \cdots, m_d} (f) \|^2_{L^p (\T^d)} \Big)^{1/2} \leq \Big(  \sum_{m_1, \cdots, m_d \in \N_0} \|\D_{m_1, \cdots, m_d} (f) \|^2_{L^{\infty} (\T^d)} \Big)^{1/2},  $$
we deduce by (\ref{Pipher}) that
\begin{equation}\label{weak_mart}
 \| f \|_{L^p (\T^d)} \leq C_d p^{d/2} \Big(  \sum_{m_1, \cdots, m_d \in \N_0} \|\D_{m_1, \cdots, m_d} (f) \|^2_{L^{\infty} (\T^d)} \Big)^{1/2}
\end{equation}
for all $p>2$. Hence, to prove that (\ref{variant_LP}) holds, it suffices, in view of (\ref{weak_mart}), to show that
$$
\Big(  \sum_{m_1, \cdots, m_d \in \N_0} \|\D_{m_1, \cdots, m_d} (f) \|^2_{L^{\infty} (\T^d)} \Big)^{1/2}  \lesssim_d \Big(  \sum_{k_1, \cdots, k_d \in \N_0} \| \widetilde{\Delta}_{k_1 , \cdots, k_d} (f) \|^2_{L^{\infty} (\T^d)} \Big)^{1/2}.
$$
This last inequality follows from the next lemma which is a $d$-dimensional analogue of \cite[Lemma 2.3]{Seeger_Trebels}.
 
\begin{lemma}\label{main_lemma} Let $\delta$ be a Schwartz function that is even, supported in $ (-4,4)$ and such that $\delta|_{[-2,2]} \equiv 1$. 

Define $\psi (\xi) = \delta(  \xi) - \delta (2 \xi)$. For $k \in \N$, put $\psi_k (\xi) = \psi (2^{-k} \xi)$ and for $k=0$, put $\psi_0 = \delta$. Consider the operator
$$ \Psi_k (f) (x) = \sum_{r \in \Z} \psi_k (r)  \widehat{f}( r ) e^{i 2 \pi r x} $$
acting on functions defined over the torus.  For $k_1 ,\cdots, k_d \in \N_0 $ we use the notation $\Psi_{k_1, \cdots, k_d} = \Psi_{k_1} \otimes \cdots \otimes \Psi_{k_d}$. 

There exists a constant $C_d>0$, depending only on the dimension $d$ and on $\psi$, such that for all $d$-tuples of non-negative integers $(m_1, \cdots, m_d)$ and $(k_1, \cdots, k_d)$ one has 
\begin{equation}\label{first}
\| \E_{m_1, \cdots, m_d}  \Psi_{k_1, \cdots, k_d}   \|_{L^{\infty} (\T^d) \rightarrow L^{\infty} (\T^d) } \leq C_d \prod_{j \in A} 2^{m_j -k_j},
\end{equation}
where $A = \big\{ j \in \{1, \cdots, d\} : m_j < k_j \big\}$  and
\begin{equation}\label{second}
\| \D_{m_1, \cdots, m_d} \Psi_{k_1, \cdots, k_d}  \|_{L^{\infty} (\T^d) \rightarrow L^{\infty} (\T^d) } \leq C_d  \prod_{j=1}^d 2^{-|k_j-m_j|} 
\end{equation}
In $(\ref{first})$ we make the convention that if $A = \emptyset$, then  $\prod_{j \in A} 2^{m_j -k_j} =1 $.
\end{lemma}

The proof of Lemma \ref{main_lemma} will be given in the next subsection. By using the above lemma and in particular estimate (\ref{second}) one can easily complete the proof of Proposition \ref{LP_product}. Towards this aim, we argue as in the proof of \cite[Proposition 2.2]{Seeger_Trebels}. More precisely, we consider a trigonometric polynomial $f $ on $\T^d$ and write $f = \sum_{k_1, \cdots, k_d \in \N_0} \widetilde{\Delta}_{k_1, \cdots, k_d} (f)$. For fixed $\eta$ (and $\phi$), if $\psi $ is as in the statement of Lemma \ref{main_lemma}, then $\psi \phi =\phi$ and hence, $\Psi_{k_1 ,\cdots, k_d}    \widetilde{\Delta}_{k_1, \cdots , k_d}  = \widetilde{\Delta}_{k_1, \cdots, k_d}$. So, by using (\ref{second}),  we obtain
\begin{align*}
\| \D_{m_1, \cdots, m_d} (f) \|_{L^{\infty} (\T^d)} & \leq   \sum_{k_1 , \cdots, k_d \in \N_0}  \| \D_{m_1, \cdots, m_d} [\widetilde{\Delta}_{m_1, \cdots, m_d} (f) ] \|_{L^{\infty} (\T^d)} \\
 &\leq \sum_{k_1, \cdots, k_d \in \N_0} \| \D_{m_1, \cdots, m_d} \Psi_{k_1 ,\cdots, k_d} \|_{L^{\infty}(\T^d) \rightarrow L^{\infty}(\T^d)} \| \widetilde{\Delta}_{k_1, \cdots, k_d} ( f ) \|_{L^{\infty} (\T^d)} \\
&\lesssim_d \sum_{k_1, \cdots, k_d \in \N_0} \big( \prod_{j=1}^d 2^{-|m_j-k_j|} \big) \| \widetilde{\Delta}_{k_1, \cdots , k_d} (f )\|_{L^{\infty} (\T^d)}
\end{align*} 
and it thus follows that
\begin{align*} 
&\Big(  \sum_{m_1, \cdots, m_d \in \N_0} \|\D_{m_1, \cdots , m_d} (f) \|^2_{L^{\infty} (\T^d)} \Big)^{1/2} \lesssim_d \\ &\Big[  \sum_{m_1 , \cdots , m_d \in \N_0} \Big( \sum_{k_1, \cdots, k_d \in \N_0} \big( \prod_{j=1}^d 2^{-|m_j-k_j|} \big) \| \widetilde{\Delta}_{k_1, \cdots , k_d} (f )\|_{L^{\infty} (\T^d)} \Big)^2 \Big]^{1/2},
\end{align*}
where the implied constant  depends only on the dimension $d$. Hence, by Minkowski's integral inequality, 
\begin{align*}
&\Big[  \sum_{m_1 , \cdots , m_d \in \N_0} \Big( \sum_{k_1, \cdots, k_d \in \N_0} \big( \prod_{j=1}^d 2^{-|m_j-k_j|} \big) \| \widetilde{\Delta}_{k_1, \cdots , k_d} (f )\|_{L^{\infty} (\T^d)} \Big)^2 \Big]^{1/2} \leq\\
& \sum_{m_1, \cdots ,m_d \in \Z} \big( \prod_{j=1}^d 2^{-|m_j|} \big)   \Big( \sum_{k_1 \geq -m_1}  \cdots \sum_{k_d \geq -m_d}  \| \widetilde{\Delta}_{k_1+m_1, \cdots, m_d+k_d} (f) \|_{L^{\infty} (\T^d)}^2 \Big)^{1/2}.
\end{align*}
Since we have
\begin{align*} & \sum_{m_1, \cdots ,m_d \in \Z} \big( \prod_{j=1}^d 2^{-|m_j|} \big)   \Big( \sum_{k_1 \geq -m_1}  \cdots \sum_{k_d \geq -m_d}  \| \widetilde{\Delta}_{k_1+m_1, \cdots, m_d+k_d} (f) \|_{L^{\infty} (\T^d)}^2 \Big)^{1/2}
\lesssim \\
&\Big( \sum_{k_1, \cdots, k_d \in \N_0}   \| \widetilde{\Delta}_{k_1, \cdots, k_d} ( f )  \|_{L^{\infty} (\T^d)}^2  \Big)^{1/2}, 
\end{align*}
the proof of Proposition \ref{LP_product} will be complete once we prove Lemma \ref{main_lemma}. This will be done in the following subsection.

\subsection{Proof of Lemma \ref{main_lemma}} The proof of this Lemma is a straightforward adaptation of \cite[Lemma 2.3]{Seeger_Trebels} to the multi-parameter setting. For the sake of simplicity, we shall only present the proof of the two-dimensional case. A similar argument establishes the higher-dimensional case.

Let $\psi$ be as in the statement of Lemma  \ref{main_lemma}. Following \cite{Seeger_Trebels}, we use the notation $\psi^{(s)} (\xi) = (i 2\pi \xi)^s \psi (\xi) $, $s \in \{-1,0,1\}$ and for $k \in \N_0$ we put
$$ \Psi^{(s)}_k (f) (x) = \sum_{r \in \Z} \psi^{(s)}(2^{-k} r) \widehat{f}( r ) e^{i 2 \pi r x}.$$
For $s=0$ we write $\psi^{(0)}=\psi$ and $\Psi^{(0)}_k = \Psi_k$. Notice that we may write 
$$ \Psi^{(s)}_k (f) (x) = K^{(s)}_k \ast f (x),$$
where $K^{(s)}_k  (x) = \sum_{r \in \Z} \psi^{(s)} (2^{-k} r )e^{ i 2 \pi r x }$. Our assumption on the support of $\psi$ implies that $K^{(s)}_k $ is in fact a trigonometric polynomial on $\T$. 
By using the Poisson summation formula, see, e.g., Corollary 2.6 in Chapter VII of \cite{SteinWeiss}, it is straightforward to see that $\| K^{(s)}_k \|_{L^1 (\T)} \lesssim_{\psi} 1$. Therefore, it follows that 
$$\| \Psi^{(s_1)}_{k_1} \otimes \Psi^{(s_2)}_{k_2}  \|_{L^{\infty}(\T^2) \rightarrow L^{\infty}(\T^2)}  = \| K^{(s_1)}_{k_1} \|_{L^1 (\T)}  \| K^{(s_2)}_{k_2} \|_{L^1 (\T)} \lesssim_{\psi} 1$$ 
for all $s_1, s_2 \in \{-1,0,1\}$ and we thus deduce that
\begin{equation}\label{standard_cond}
 \sum_{s_1, s_2  \in \{-1,0,1\} } \| \Psi^{(s_1)}_{k_1} \otimes  \Psi^{(s_2)}_{k_2} ( f) \|_{L^{\infty}(\T^2)} \lesssim_{\psi} \| f \|_{L^{\infty}(\T^2)} \\
 \end{equation}
 for all $k_1,k_2 \in \N_0$, where the summation is taken with respect to all possible choices of $s_1, s_2 \in \{-1,0,1\}$.
 \subsubsection{Proof of condition $(\ref{first})$ (for $n=2$).} 
We shall consider two cases; $A =\emptyset$ and $A \neq \emptyset$.\newline
\underline{Case 1: $A = \emptyset$.} In this case we have $m_1 \geq k_1$ and $m_2 \geq k_2$ and (\ref{first}) easily follows from (\ref{standard_cond}), 
$$ \|\E_{m_1, m_2} \Psi_{k_1, k_2} \|_{L^{\infty}(\T^2) \rightarrow L^{\infty}(\T^2)} \lesssim 1.$$
\underline{Case 2: $A \neq \emptyset$.}
First, consider the subcase where $m_1<k_1$ and $m_2< k_2 $. For $(x_1, x_2) \in \T^2$, we denote by $I_j = [s_j 2^{-m_j}, (s_j + 1) 2^{-m_j} )$, $s_j \in \{0,1 ,\cdots, 2^{m_j} -1 \}$, the unique dyadic interval in $\T$ of length $2^{-m_j}$ containing $x_j$  ($j=1,2$). If we write $I_j = [a_j, b_j)$, i.e. $a_j = s_j 2^{-m_j}$, $b_j = (s_j+1) 2^{-m_j}$, then we have
\begin{align*}
& \E_{m_1, m_2} [ \Psi_{k_1, k_2} ( f ) ] (x_1, x_2) = \\
& 2^{m_1} 2^{m_2} \int_{I_1 \times I_2}  \sum_{r_1 , r_2 \in \Z} \psi (2^{-k_1} r_1) \psi(2^{-k_2} r_2) \widehat{f} (r_1,r_2) e^{i2\pi (r_1 x'_1 + r_2 x'_2)} d x'_1 d x'_2 =\\& 2^{m_1} 2^{m_2}  \sum_{r_1 , r_2 \in \Z}  \psi (2^{-k_1} r_1) \psi(2^{-k_2} r_2) \widehat{f} (r_1, r_2) 
\big[ \frac{ e^{i 2 \pi r_1 b_1}   - e^{i 2 \pi r_1 a_1}  } {i 2 \pi r_1} \big] \big[  \frac{ e^{i 2 \pi r_2 b_2} - e^{i 2 \pi r_2 a_2 } } {i 2 \pi r_2} \big] .
\end{align*}
Hence, one can write
\begin{align*}
& \E_{m_1, m_2} [ \Psi_{k_1, k_2}  (f)] (x_1, x_2) =\\
& 2^{-k_1+m_1} \cdot 2^{-k_2+m_2} \big[ \Psi_{k_1}^{(-1)} \otimes \Psi_{k_2}^{(-1)} ( f ) (b_1,b_2) - \Psi_{k_1}^{(-1)} \otimes \Psi_{k_2}^{(-1)} ( f) (b_1, a_2) \\
&- \Psi_{k_1}^{(-1)} \otimes \Psi_{k_2}^{(-1)} (f) (a_1, b_2) + \Psi_{k_1}^{(-1)} \otimes \Psi_{k_2}^{(-1)} ( f ) (a_1, a_2) \big] 
\end{align*}
and thus, by (\ref{standard_cond}), we obtain the desired estimate,
$$ \|\E_{m_1, m_2} \Psi_{k_1 , k_2} \|_{L^{\infty}(\T^2) \rightarrow L^{\infty}(\T^2)} \lesssim 2^{m_1- k_1} \cdot 2^{ m_2 -k_2}. $$
Next, consider the subcase where $m_1<k_1$ but $m_2 \geq k_2$. In this case, for $(x_1,x_2) \in I_1 \times I_2 = [a_1,b_1) \times [a_2, b_2)$, $I_1$, $I_2$ being as in the previous subcase, we have
\begin{align*}
& \E_{m_1, m_2} [ \Psi_{k_1, k_2} ( f ) ] (x_1, x_2)  =\\
& 2^{m_1 -k_1} 2^{m_2} \int_{I_2} \big[  \sum_{r_1 , r_2 \in \Z}  \psi (2^{-k_1} r_1) \psi(2^{-k_2} r_2) \widehat{f} (r_1, r_2)  e^{i 2 \pi r_2 x'_2} \big( \frac{e^{i 2 \pi r_1  b_1}  - e^{i 2 \pi r_1 a_1 } } {i 2 \pi r_1} \big)  \big] dx'_2 ,
\end{align*}
and so, $\E_{m_1, m_2} [\Psi_{k_1, k_2} (f)] (x_1,x_2)$ can be written as
$$ 2^{m_1-k_1}  2^{m_2} \big[ \int_{ I_2 }  \Psi^{(-1)}_{k_1} \otimes \Psi_{k_2} ( f)  (b_1,x'_2) - \Psi^{(-1)}_{k_1} \otimes \Psi_{k_2} ( f ) (a_1, x'_2) d x'_2 \big].  $$
Since the length of $I_2$ is equal to $2^{-m_2}$, we get
$$ \|  \E_{m_1, m_2} [\Psi_{k_1, k_2} ( f) ] \|_{L^{\infty}(\T^2)} \leq 2^{m_1-k_1}  \cdot 2^{m_2} \cdot  2^{-m_2} \cdot  2 \|  \Psi^{(-1)}_{k_1} \otimes \Psi_{k_2} (f)   \|_{L^{\infty} (\T^2)} $$
and hence, by using (\ref{standard_cond}), we have 
$$ \|\E_{m_1, m_2} \Psi_{k_1, k_2} \|_{L^{\infty}(\T^2) \rightarrow L^{\infty}(\T^2)} \lesssim 2^{m_1 -k_1}  . $$
The subcase where $m_1 \geq k_1$ and $m_2< k_2$ is symmetric to the previous one. Therefore, (\ref{first}) is completely shown in the two-dimensional case. \newline

\subsubsection{Proof of condition $(\ref{second})$ (for $n=2$).} 
We shall consider two cases; $A = \{1,2\}$ and $ \{1,2\} \setminus A \neq \emptyset$.\newline
 \underline{Case 1: $A = \{1,2\}$.} In this case we have $m_1 < k_1$ and $m_2 < k_2$ and (\ref{second}) follows easily from (\ref{first}). Indeed, observe that
\begin{align*}
\| \D_{m_1, m_2} [ \Psi_{k_1, k_2} (f) ]\|_{L^{\infty}(\T^2)}&\leq \|\E_{m_1,m_2} [\Psi_{k_1, k_2} (f) ] \|_{L^{\infty}(\T^2)} + \|\E_{m_1-1,m_2}  [ \Psi_{k_1, k_2}  (f) ] \|_{L^{\infty}(\T^2)} \\
& +\|\E_{m_1,m_2-1} [ \Psi_{k_1, k_2}  (f) ]  \|_{L^{\infty}(\T^2)} + \|\E_{m_1-1,m_2-1} [ \Psi_{k_1, k_2}  (f) ] \|_{L^{\infty}(\T^2)} \\
&\lesssim 2^{m_1-k_1} \cdot 2^{m_2- k_2} \| f \|_{L^{\infty}(\T^2)}
\end{align*}
by (\ref{first}), as $m_1-1<m_1<k_1$ and $m_2-1<m_2<k_2$. \newline
\underline{Case 2: $\{1,2\} \setminus A \neq \emptyset$.}
Assume first that $A = \emptyset$, that is $m_1 \geq k_1$ and $m_2 \geq k_2$. By using the definition of $\D_{m_1, m_2}$, we write
\begin{align*}
 \D_{m_1,m_2} [ \Psi_{k_1, k_2} (f) ] &= \E_{m_1, m_2}  [ \Psi_{k_1, k_2} ( f) ] -  \E_{m_1, m_2 -1}  [ \Psi_{k_1, k_2} ( f) ]\\
 &- ( \E_{m_1 -1, m_2 }  [ \Psi_{k_1, k_2} ( f) ] - \E_{m_1 -1 , m_2 -1}  [ \Psi_{k_1, k_2} ( f) ] ). 
 \end{align*}
Take $(x_1, x_2) \in \T^2$ and for $j=1,2$ let $I_j$ be the dyadic interval in $\T$ of length $2^{-m_j}$ containing $x_j$. Let $\widetilde{I}_j$ denote the dyadic interval of length $2^{-m_j +1}$ such that $x_j \in \widetilde{I}_j$. Note that since $I_j$ and $\widetilde{I}_j$ are dyadic intervals with non-empty intersection and $|\widetilde{I}_j| = 2 |I_j|$, one has $I_j \subset \widetilde{I}_j$, $j=1,2$. Since
\begin{align*}
&\E_{m_1, m_2}  [ \Psi_{k_1, k_2} ( f) ] (x_1, x_2) -  \E_{m_1 , m_2 - 1}  [ \Psi_{k_1, k_2} ( f) ] (x_1, x_2) =\\
&2^{-m_1} \int_{I_1} \big( \E_{m_2} \Psi_{k_1, k_2} (f) (x'_1, x_2) -  \E_{m_2 - 1} \Psi_{k_1, k_2} (f) (x'_1, x_2) \big) dx'_1,
\end{align*}
by using the mean value theorem for integrals it follows that there exists an $x^{(\alpha)}_1 \in I_1$ such that
\begin{align*}  \E_{m_1, m_2}  [ \Psi_{k_1, k_2} ( f) ] (x_1, x_2) -  &\E_{m_1 , m_2 - 1}  [ \Psi_{k_1, k_2} ( f) ] (x_1, x_2)  =\\
&  \E_{m_2} \Psi_{k_1, k_2} (f) (x^{(\alpha)}_1 , x_2) -  \E_{m_2 - 1} \Psi_{k_1, k_2} (f) (x^{(\alpha)}_1 , x_2).
\end{align*}
A similar analysis on $\E_{m_1 -1, m_2}  [ \Psi_{k_1, k_2} ( f) ] -  \E_{m_1 -1 , m_2 - 1}  [ \Psi_{k_1, k_2} ( f) ] $ shows that there exists an $x^{(\beta)}_1 \in \widetilde{I}_1$ such that
\begin{align*}  \E_{m_1 -1, m_2}  [ \Psi_{k_1, k_2} ( f) ] (x_1, x_2) -  &\E_{m_1 -1 , m_2 - 1}  [ \Psi_{k_1, k_2} ( f) ] (x_1, x_2)  =\\
&  \E_{m_2} \Psi_{k_1, k_2} (f) (x^{(\beta)}_1 , x_2) -  \E_{m_2 - 1} \Psi_{k_1, k_2} (f) (x^{(\beta)}_1 , x_2).
\end{align*}
Therefore,
\begin{align*}
\D_{m_1, m_2} [\Psi_{k_1, k_2} (f)] (x_1, x_2) &=  \E_{m_2} \Psi_{k_1, k_2} (f) (x^{(\alpha)}_1 , x_2) -  \E_{m_2 - 1} \Psi_{k_1, k_2} (f) (x^{(\alpha)}_1 , x_2) \\
& - \big(  \E_{m_2} \Psi_{k_1, k_2} (f) (x^{(\beta)}_1 , x_2) -  \E_{m_2 - 1} \Psi_{k_1, k_2} (f) (x^{(\beta)}_1 , x_2) \big) . 
\end{align*}
If we assume, without loss of generality, that $x^{(\alpha)}_1 < x^{(\beta)}_1$, then by the mean value theorem, 
$$ \D_{m_1, m_2} [\Psi_{k_1, k_2} (f)] (x_1, x_2) =  (x^{(\beta)}_1 - x^{(\alpha)}_1 )\partial_{x_1} \big\{ \E_{m_2 } \Psi_{k_1, k_2} (f)  - \E_{m_2 -1} \Psi_{k_1, k_2} (f)  \big\} (x^{(\gamma)}_1, x_2)
$$for some $x^{(\gamma)}_1 \in (x^{(\alpha)}_1 , x^{(\beta)}_1)$. 
One can easily see that
$$\partial_{x_1} \big\{ \E_{m_2 } \Psi_{k_1, k_2} (f)  - \E_{m_2 -1} \Psi_{k_1, k_2} (f)  \big\}  = 2^{k_1} \big( \E_{m_2 } \Psi^{(1)}_{k_1} \otimes \Psi_{k_2} (f)  -  \E_{m_2 - 1} \Psi^{(1)}_{k_1} \otimes \Psi_{k_2} (f) \big) $$
and so,
$$ \D_{m_1, m_2} [\Psi_{k_1, k_2} (f)] (x_1, x_2) =   2^{k_1} \big( \E_{m_2 } \Psi^{(1)}_{k_1} \otimes \Psi_{k_2} (f)  (x^{(\gamma)}_1, x_2)  -  \E_{m_2 - 1} \Psi^{(1)}_{k_1} \otimes \Psi_{k_2} (f) (x^{(\gamma)}_1, x_2) \big).
$$
Hence, by using the definition of $E_{m_2}$ and $E_{m_2 - 1}$, it follows by the mean value theorem for integrals that there are $x^{(\alpha)}_2 \in I_2$ and $x^{(\beta)}_2 \in \widetilde{I}_2$ such that
\begin{align*}
& \D_{m_1, m_2} [\Psi_{k_1, k_2} (f)] (x_1, x_2) =   \\
& (x^{(\beta)}_1 - x^{(\alpha)}_1 ) 2^{k_1} \big( \Psi^{(1)}_{k_1} \otimes \Psi_{k_2} (f)  (x^{(\gamma)}_1, x^{(\alpha)}_2)  -   \Psi^{(1)}_{k_1} \otimes \Psi_{k_2} (f) (x^{(\gamma)}_1, x^{(\beta)}_2) \big).
\end{align*}
Without loss of generality we may assume that $x^{(\alpha)}_2 < x^{(\beta)}_2$. Hence, by applying the mean value theorem, we deduce that
$$  \D_{m_1, m_2} [\Psi_{k_1, k_2} (f)] (x_1, x_2) =   (x^{(\beta)}_1 - x^{(\alpha)}_1 ) (x^{(\beta)}_2 - x^{(\alpha)}_2) 2^{k_1} 2^{k_2}  \Psi^{(1)}_{k_1} \otimes \Psi^{(1)}_{k_2} (f)  (x^{(\gamma)}_1, x^{(\gamma)}_2) 
$$
for some $x^{(\gamma)}_2 \in ( x^{(\alpha)}_1, x^{(\beta)}_2)$. Since $|x^{(\beta)}_j - x^{(\alpha)}_j | \leq 2^{-m_j +1}$, we obtain
$$  \| \D_{m_1, m_2} [\Psi_{k_1, k_2} (f)]  \|_{L^{\infty}(\T^2) \rightarrow L^{\infty} (\T^2) }  \lesssim 2^{k_1- m_1} 2^{k_2 -m_2} \| f \|_{L^{\infty} (\T^2)}, $$
as desired.

It only remains to consider the subcase where $m_1 \geq k_1$ and $m_2 < k_2$, the other one ($m_1 < k_1$ and $m_2 \geq k_2$) being symmetric. We need to show that
$$ \| \D_{m_1, m_2} \Psi_{k_1, k_2} \|_{L^{\infty}(\T^2) \rightarrow L^{\infty} (\T^2)} \lesssim 2^{k_1 - m_1}  2^{ m_2 - k_2 }. $$ 
To this end, write $ \D_{m_1,m_2} [ \Psi_{k_1, k_2} (f) ] $ as $$ \E_{m_1, m_2}  [ \Psi_{k_1, k_2} ( f) ] -  \E_{m_1 -1, m_2 }  [ \Psi_{k_1, k_2} ( f) ]
 - ( \E_{m_1, m_2 - 1 }  [ \Psi_{k_1, k_2} ( f) ] - \E_{m_1 -1 , m_2 -1}  [ \Psi_{k_1, k_2} ( f) ] ) $$
 and handle each of these two terms separately. Take $(x_1, x_2) \in \T^2$ and, for $j=1,2$, consider the dyadic intervals $I_j = [a_j, b_j)$ and $\widetilde{I}_i$ as above.  For the first term, by applying the mean value theorem for integrals, we see that there are $x^{(\alpha)}_1 \in I_1$  and $x^{(\beta)}_1 \in \widetilde{I}_1$ such that
 \begin{align*}
  \E_{m_1, m_2}  [ \Psi_{k_1, k_2} ( f) ] (x_1, x_2) -  \E_{m_1 -1, m_2 }  [ \Psi_{k_1, k_2} ( f) ] (x_1, x_2) &=  \E_{m_2}   \Psi_{k_1, k_2} ( f)  (x^{(\alpha)}_1, x_2) \\
  &-  \E_{ m_2 }   \Psi_{k_1, k_2} ( f)  (x^{(\beta)}_1, x_2) . 
 \end{align*}
 Hence, if we assume that $x^{(\alpha)}_1 <x^{(\beta)}_1$, then by the mean value theorem there is an $x^{(\gamma)}_1 \in ( x^{(\alpha)}_1, x^{(\beta)}_1)$ such that
  \begin{align*}
  \E_{m_1, m_2}  [ \Psi_{k_1, k_2} ( f) ] (x_1, x_2) &-  \E_{m_1 -1, m_2 }  [ \Psi_{k_1, k_2} ( f) ] (x_1, x_2)   \\
  & = ( x^{(\beta)}_1  -  x^{(\alpha)}_1 ) \partial_{x_1} \E_{ m_2}   \Psi_{k_1, k_2} ( f)  (x^{(\gamma)}_1 , x_2) \\
  &= ( x^{(\beta)}_1  -  x^{(\alpha)}_1 ) 2^{k_1} \E_{ m_2}   \Psi^{(1)}_{k_1} \otimes \Psi_{ k_2} ( f)  (x^{(\gamma)}_1 , x_2).
 \end{align*}
 Now, by considering the definition of $\E_{m_2}$, an explicit calculation shows that
 \begin{align*}
  &\E_{m_1, m_2}  [ \Psi_{k_1, k_2} ( f) ] (x_1, x_2) -  \E_{m_1 -1, m_2 }  [ \Psi_{k_1, k_2} ( f) ] (x_1, x_2)   \\
  &= ( x^{(\beta)}_1  -  x^{(\alpha)}_1 ) 2^{k_1}  2^{m_2} 2^{-k_2} \big[   \Psi^{(1)}_{k_1} \otimes \Psi^{(-1)}_{ k_2} ( f)  (x^{(\gamma)}_1 , b_2) -  \Psi^{(1)}_{k_1} \otimes \Psi^{(-1)}_{ k_2} ( f)  (x^{(\gamma)}_1 , a_2) \big],
 \end{align*}
 where $I_2 = [a_2, b_2)$ is as above. Since $ x^{(\alpha)}_1, x^{(\beta)}_1 \in \widetilde{I}_1$, the last expression gives
 $$ \| \E_{m_1, m_2}  [ \Psi_{k_1, k_2} ( f) ]  -  \E_{m_1 -1, m_2 }  [ \Psi_{k_1, k_2} ( f) ] \|_{L^{\infty} (\T^2)} \lesssim    2^{k_1 - m_1} 2^{m_2-k_2} \| f \|_{L^{\infty} (\T^2)}.$$
A similar argument shows that the second term also satisfies
$$ \| \E_{m_1, m_2 - 1 }  [ \Psi_{k_1, k_2} ( f) ] - \E_{m_1 -1 , m_2 -1}  [ \Psi_{k_1, k_2} ( f) ] \|_{L^{\infty} (\T^2)} \lesssim    2^{k_1 - m_1} 2^{m_2-k_2} \| f \|_{L^{\infty} (\T^2)} $$
and we thus deduce that
$$ \| \D_{m_1, m_2} \Psi_{k_1, k_2} \|_{L^{\infty}(\T^2) \rightarrow L^{\infty} (\T^2)} \lesssim 2^{k_1 - m_1}  2^{ m_2 - k_2 } . $$
Hence, the proof of (\ref{second}) for $n=2$ is complete.

\section{Some Further remarks and applications}  
\subsection{Applications in thin sets}

Proposition \ref{2d_Zygmund} gives examples of $\Lambda(p)$ sets in $\Z^d$ whose corresponding $\Lambda(p)$ constant grows like $p^{d/2}$ as $p \rightarrow \infty$ and they cannot be written as products of Sidon sets. Moreover, those sets, namely the class of the sets $E \subset \Z^d$ that cannot be written as $d$-fold products of sets in $\Z$ and satisfy the condition $\sup_{I_1 , \cdots, I_d \in \mathcal{J}} \# \big\{ E \cap (I_1 \times \cdots \times I_d) \big\} < \infty $, are examples of $2d/(d+1)$-Rider sets in $\Z^d$ that cannot be written as  products of Sidon sets in $\Z$.

Note that if $\Lambda_1,  \cdots, \Lambda_d$ are lacunary sequences in $\Z$, then $ \Lambda_1 \times \cdots \times \Lambda_d$ satisfies (\ref{rectangle}) and we thus recover the well-known fact that $ \Lambda_1 \times \cdots \times \Lambda_d$ is a $\Lambda(p)$ set in $\Z^d$ whose constant grows like $p^{d/2}$ as $p \rightarrow \infty$. However,  Proposition \ref{2d_Zygmund} cannot handle spectral sets of the form $\Lambda_1 \times \cdots \times \Lambda_d$, where $\Lambda_j$ is a Sidon set that is not a finite union of lacunary sequences ($j=1, \cdots, d$).

\subsection{A version of (\ref{variant_LP}) for ``rough'' projections} For $k \in \N$ consider the classical Littlewood-Paley projections
$$ \Delta_k (f)(x) = \sum_{n=2^{k - 1}}^{2^k -1} \widehat{f} (n) e^{ i 2 \pi n x} +  \sum_{n=-2^k +1 }^{-2^{k-1}} \widehat{f} (n) e^{ i 2 \pi n x}   .$$
For $k =0$, set $ \Delta_0 (f) (x) = \widehat{f} (0)$. For $k_1, \cdots, k_d \in \N_0$ we write 
$$ \Delta_{k_1, \cdots ,k_d} =  \Delta_{k_1} \otimes \cdots \otimes  \Delta_{k_d} . $$
Since for every trigonometric polynomial $f$ on the $d$-torus we may write $f = \sum_{m_1, \cdots, m_d \in \N_0}  \Delta_{m_1, \cdots, m_d} (f)$, 
we have
$$ \widetilde{ \Delta}_{k_1, \cdots, k_d} (f) = \sum_{m_1, \cdots, m_d \in \N_0} \widetilde{ \Delta}_{k_1, \cdots, k_d}  \Delta_{m_1, \cdots, m_d} (f).$$
Observe that $\widetilde{ \Delta}_{k_1, \cdots, k_d}  \Delta_{m_1, \cdots, m_d} = 0$ whenever there exists an index $j_0 \in \{1, \cdots, d\}$ such that $|k_{j_0}- m_{j_0}  |>1$. We thus deduce that
\begin{align*} \| \widetilde{ \Delta}_{k_1, \cdots, k_d} (f) \|_{L^{\infty} (\T^d)} &\leq \sum_{\substack{(m_1, \cdots, m_d) \in \N_0^d: \\ |k_j- m_j | \leq 1 \ \mathrm{for} \ \mathrm{all} \ j \in \{1, \cdots, d\} } } \| \widetilde{ \Delta}_{k_1, \cdots, k_d}  \Delta_{m_1, \cdots, m_d} (f) \|_{L^{\infty} (\T^d)} \\
&\lesssim_d \sum_{\substack{(m_1, \cdots, m_d) \in \N_0^d: \\ |k_j- m_j | \leq 1 \ \mathrm{for} \ \mathrm{all} \ j \in \{1, \cdots, d\} } }  \|   \Delta_{m_1, \cdots, m_d} (f) \|_{L^{\infty} (\T^d)}.
\end{align*}
Therefore,
$$ \Big( \sum_{k_1, \cdots, k_d \in \N_0} \| \widetilde{ \Delta}_{k_1, \cdots, k_d} (f) \|^2_{L^{\infty} (\T^d)} \Big)^{1/2} \lesssim_d \Big( \sum_{k_1, \cdots, k_d \in \N_0} \|  \Delta_{k_1, \cdots, k_d } (f) \|^2_{L^{\infty} (\T^d)} \Big)^{1/2} $$
and hence, it follows by (\ref{variant_LP}) that for every trigonometric polynomial $f$ on $\T^d$ one has
\begin{equation}\label{rough_LP}
\| f \|_{L^p (\T^d)} \lesssim_d p^{d/2} \Big( \sum_{k_1, \cdots, k_d \in \N_0} \|  \Delta_{k_1, \cdots, k_d} (f) \|^2_{L^{\infty} (\T^d)} \Big)^{1/2}
\end{equation}
for every $p>2$. Estimate (\ref{rough_LP}) is a multi-parameter version of an inequality due to C. Moore \cite{Moore_thesis}. In particular, we obtain the following multi-parameter extension of \cite[Theorem, p.30]{Moore_thesis}.

\begin{corollary}
There exist positive constants $c_1 (d)$ and $c_2 (d)$, depending only on the dimension $d$, such that whenever 
$$  \sum_{k_1, \cdots, k_d \in \N_0 } \|  \Delta_{k_1, \cdots, k_d} (f) \|^2_{L^{\infty} (\T^d)} < \infty$$
one has
$$ \int_{\T^d} \exp \Big\{ c_1 (d) \Big[ \frac{|f(x_1, \cdots, x_d)|}{\big( \sum_{ k_1, \cdots, k_d \in \N_0 } \|  \Delta_{k_1, \cdots, k_d} (f) \|^2_{L^{\infty} (\T^d)} \big)^{1/2} } \Big]^{2/d} \Big\} dx_1 \cdots d x_d < c_2 (d).$$
\end{corollary}

\bibliographystyle{plainnat}
\bibliography{m-i_biblio}

\end{document}